\input amstex
\documentstyle{amsppt}
\magnification=1200
\overfullrule=0pt
\NoRunningHeads
%
\def\C{\Bbb C}

\def\H{\Bbb H}
\def\K{\Bbb K}

\def\P{\Bbb P}
\def\R{\Bbb R}

\def\Z{\Bbb Z}

%
\topmatter
\title
Euler class and free generation
\endtitle
%
\author
Alexander Reznikov 
\endauthor
%
%
%
%
\address
Institute of Mathematics, Hebrew University, Giv'at Ram 91904, Jerusalem, Israel.
\endaddress
%
%
%
%
\date
July, 1996
\enddate
%
%
%
%
\thanks
Partially supported by a BSF grant
\endthanks
%
%
%
%
%
%
%
%
%
%

\endtopmatter
\document

$   $\newline
``What can you say in your favour?'' \newline
``You see ...'' \newline
``Enough. Shoot'm. Next.'' [S].

$   $\newline

This paper consists of two parts. In the first auxilliary part, we deal with sets with cyclic order. For any such set
$\Cal O$
we introduce following [BG] a cocycle
$\ell : G \times G \rightarrow \Z$
valued in
$\{0, \pm 1\} \subset \Z$
on the group $G$ of automorphisms of
$\Cal O$.
The cohomology class of
$\ell$
in
$H^2(G,\Z)$
will be called the Euler class. If $K$ is an ordered field, then the projective line
$\P^1(K)$
has a cyclic order and
$PSL_2(K)$
acts order-preserving on
$\P^1(K)$, so that we get both the cocycle 
$\ell$
and the Euler class in
$H^2(PSL_2(K),\Z)$.
If
$K = \R$,
then our Euler class coincides with the usual Euler class on
$PSL_2(\R)$.

In view of our extension of the Euler class to all ordered fields, the following two problems arise.

\demo{{\bf Problem 1}}
Let
$\rho : \pi_1(S) \rightarrow PSL_2(K)$
be a homomorphism of the fundamental group of a closed oriented surface. Is it true, that
$|(\rho^*[\ell], [S])| \le 2g-2$?
\enddemo

\demo{{\bf Problem 2}}
Suppose
$|(\rho^*[\ell], [S])| = 2g-2$.
Is it true that
$\rho$
is injective?
\enddemo

For
$\K = \R$
the theorems of Milnor [M] and Goldman [Go2], answer these problems positively.

In the main second part of this paper, we apply the cocycle
$\ell$
and the ideas from the theory of Hamiltonian systems on the Teichm\"uller space to the following classical problem:

\demo{{\bf Problem 3}}
When $n$ matrices in
$SL_2(\R)$
generate a free discrete group?
\enddemo

For
$n=2$
this problem has been treated in many papers, see [Gi]. An effective solution is given in [Gi]. The analogous problem for
$SL_2(\C)$
also attracted a lot of attention especially since J\"orgensen paper [J]. For
$n>2$
however, the problem becomes much harder. We will give a simple {\it sufficient} condition for $n$ hyperbolic elements in
$SL_2(\R)$
to generate a free discrete group. This condition is open, that is, satisfied on an open domain in
$(SL_2(\R))^n$.
Here is our main result.

Let $n$ be even and let
$a_i,b_i, 1\le i \le n$
be in
$SL_2(\R)$.
Suppose
$h = \prod^n_{i=1} [a_i,b_i]$
is hyperbolic. Consider the eigenvectors
$x_1, x_2$
of $h$ and a matrix $r$ which takes the form
$\pmatrix 1 & 0 \\ 0 & -1 \endpmatrix$
in the basis
$(x_1, x_2)$.
Put
$a_i = rb_{n+1-i}r^{-1}, b_i = ra_{n+1-i}r^{-1}$
for
$n+1 \le i \le 2n$.
Let
$I_j = \prod^j_{i=1} [a_i,b_i] (j \le 2n)$.

\proclaim{Main Theorem}
Let
$f(a,b) = \frac 1 {\pi} \sum^{2n}_{j=1} \ell (I_{j-1},a_j) + \ell (I_{j-1}a_j,b_j)- \ell (I_{j-1}a_jb_j a^{-1}_j,a_j) - \ell (I_j,b_j)$
Then
\roster
\item"(a)" $f(a,b)$
is an integer and
$|f(a,b)| \le 2n-1$
\item"(b)" if
$|f(a,b)| = 2n-1$,
then
$\{a_i,b_i\}$
generate a free hyperbolic group in
$SL_2(\R)$.
\endroster
\endproclaim

\demo{\bf{1. Cyclically ordered sets, ordered fields and the Euler class}} 
\enddemo

\demo{1.1}
A cyclically ordered set
$\Cal O$
is a set with a subset
$\Omega$
in
$\Cal O \times \Cal O \times \Cal O$,
satisfying the following conditions:
\roster
\item"(i)" if
$(x,y,z) \in \Cal O$
then
$x,y,z$
are all different
\item"(ii)" if
$\sigma$
is a permutation of
$(x,y,z)$
and
$(x,y,z) \in \Cal O$
then
$\{\sigma (x), \sigma (y), \sigma (z)\} \in \Cal O$
if and only if $\sigma$ is even
\item"(iii)" if $z$ is fixed then the relation
$x < y \Leftrightarrow (x,y,z) \in \Cal O$
is a linear order.
\endroster
\enddemo

\demo{1.2 Example} Let $K$ be an ordered field and let
$\P^1(K)$
be a projective line over $K$. We can think of
$\P^1 (K)$
as
$K \cup\{\infty\}$.
The cyclic order in
$\P^(K)$
is defined by a condition that the induced order in $K$ is standard. The group
$PSL_2(K)$
acts on
$\P^1(K)$
preserving the cyclic order.
\enddemo

\demo{1.3}
Define a function
$\psi : \Cal O \times \Cal O \times \Cal O \rightarrow \{0, \pm 1\}$
in a following way:
\roster
\item"(i)" if any of
$(x,y,z)$
are equal, then
$\psi (x,y,z) =0$
\item"(ii)"
$\psi$
is odd under permutation of
$(x,y,z)$
\item"(iii)" if
$(x,y,z) \in \Omega$,
then
$\psi(x,y,z) =1$.
\endroster
\enddemo

\demo{1.4} Now let $G$ be a group, acting in order preserving way on
$\Cal O$.
Fix any element
$p \in \Cal O$
and define a function
$\ell : G \times G \rightarrow \{0, \pm 1\}$
as
$\ell(g_1,g_2) = \psi (p,g_2p,g_1g_2p)$.
\enddemo

\proclaim{Lemma (1.4)}
$\ell$
is an integer cocycle on $G$.
\endproclaim

\demo{Proof}
is a direct computation and left to the reader.
\enddemo

\demo{Definition}
The cohomology class of
$\ell$
in
$H^2(G)$
(which does not depend on $p$) is called the Euler class. In particular, for an ordered field $K$ one gets the Euler class in
$H^2(PSL_2(K), \Z)$.
\enddemo

\proclaim{1.5 Comparison theorem ([BG])}
For
$K = \R$,
the class of
$\ell$
in
$H^2(PSL_2(\R))$
is the usual Euler class of associated
$S^1$-bundle over
$BPSL_2^{\delta}(\R)$.
\endproclaim

\demo{Proof}
Consider the action of
$PSL_2(\R)$
on
$\Cal H^2$.
For any
$p \in \Cal H^2$
the class
$\ell_p(g_1,g_2) = \text{Area}\, (p,g_2p,g_1g_2p)$
represents the Euler class $e$ [Gu]. Here
$\text{Area} \, (p,q,r)$
is the area of oriented geodesic triangle with vertices in
$p,q,r$.
Now
$[\ell_p] \in H^2(PSL_2(\R))$
does not depend on $p$ and all cocycles 
$\ell_p$
are uniformly bounded. For
$p_0 \in \partial \Cal H^2$
and
$p \rightarrow p_0$,
we will have
$\ell^{\infty}$-convergence of cocycles
$\ell_{p_i} \rightarrow \ell_{p_0}$.
Moreover the area of ideal triangle
$(p,q,r)$
is
$\pi \cdot \psi(p,q,r)$.
Any homology class in
$H_2(PSL_2(\R))$
is represented by a map of a surface group
$\pi_1(S) \overset{\alpha} \to \longrightarrow PSL_2(\R)$.
It follows that
$(\ell_{p_0}, \alpha_*[S]) = \lim(\ell_{p_i}, \alpha_*[S]) = (e,\alpha_*[S])$.
This completes the proof.
\enddemo

\demo{\bf {2. Discrete Goldman twist}}
We will work with the representation variety
$\Cal M = \text{Hom} \, (\pi_1(S), SL_2(\R))/SL_2(\R)$,
where $S$ is an oriented closed surface of genus $g$. It has
$2^{2g+1}+2g-3$
connected components, indexed by the value of the Euler class [H]. Every such component, say
$\Cal M_e$,
is a symplectic manifold, nonsingular if
$e \ne 0$.
For any
$\gamma$
a conjugacy class in
$\pi_1(M)$,
there is a natural Hamiltonian 
$Tr_{\gamma} : \Cal M \rightarrow \R$,
and the corresponding Hamiltonian flow has been identified by Goldman. If
$\gamma$
can be represented by a simple separating curve, this flow can be described as follows. Write a presentation of
$\pi_1(S)$
in the following form:

$$[x_1,x_2] \ldots [x_{2\kappa -1},x_{2\kappa}] = [\gamma] = [x_{2\kappa+1}, x_{2\kappa+2}] \ldots [x_{2g-1},x_{2g}]$$
Next, write
$[\gamma] = \exp A$
for some
$A \in sl_2(\R)$.
Then put

$$\aligned \bar x^t_i &= \bar x_i, i \le 2\kappa \\
\bar x^t_i &= \exp(-tA) \bar x_i \exp tA, i \ge 2\kappa; \endaligned$$
this is the flow of
$Tr_{\gamma}$.
In particular,
$f_t: \{\bar x_i\} \rightarrow \{\bar x^t_i\}$
is a symplectomorphism of
$\Cal M$.
Here
$\bar x_i$
stands for the representation matrix of
$x_i$.

For different
$\gamma$,
the Hamiltonians
$Tr_{\gamma}$
yield nice commutation relations, discovered by Wolpert [W] and put in a more ``representation variety language'' by Goldman [Go1]. In fact, the free module on the set of conjugacy classes  of
$\pi_1(M)$
becomes a Lie ring with Goldman's bracket. One may wonder what kind of group object correspond to it.

Whatever this eventual ``Kac-Moody-Goldman'' group may be, we will introduce now some elements from the ``other connected components'' of it. These are defined for
$\Cal M_{\pm (g-1)}$,
which is naturally symplectomorphic to the Teichm\"uller space by a theorem of Goldman (see various proofs in [Go2], [H], [Re2]. In this case, all representation matrices are hyperbolic.

For
$\gamma$
as above, let $r$ be a unique matrix (up to sign), commuting with
$[\bar {\gamma}]$
with eigenvalues
$+1$
and
$-1$.
Put

$$\aligned f(\bar x_i) &= x_i, i \le 2\kappa \\
f(\bar x_i) &= r^{-1} \bar x_i r, i \ge 2\kappa \endaligned$$
This is a symplectic diffeomorphism of
$\Cal M_{\pm(g-1)}$.
There is a particularly nice description of the map $f$ if one views
$\Cal M_{g-1}$
as Teichm\"uller space. Namely, realise a point in
$\Cal M_{g-1}$
as a hyperbolic metric on $S$ and find a geodesic, representing
$\gamma$.
We assume that the marked point lies on
$\gamma$.
Cut $S$ into two pieces along
$\gamma$
and glue again by a reflection, which fixes a marked point. Then the new hyperbolic structure is the image of $f$.
\enddemo

\demo{\bf{3. Euler class}}
For a representation
$x_i \rightarrow \bar x_i$
of
$\pi_1(S)$
in
$SL_2(\R)$
the Euler number is an integer between
$-(g-1)$
and
$(g-1)$
by Milnor [M]. As mentioned above, all representations with the maximal Euler number are discrete faithful hyperbolic by Goldman's theorem. One can introduce a universal Euler class $e$ in
$H^2(SL^{\delta}_2(\R),\R)$
as the image of a generator in continuous cohomology
$H^2_{cont}(SL_2(\R),\R)$.
A representation
$x_i \rightarrow \bar x_i$
as above defines a homology class in
$H_2(SL^{\delta}_2(\R),\Z)$,
the image of the generator of
$H_2(\pi_1(S),\Z) \approx \Z$,
and the Euler number is just given by the evaluation of $e$ on this class. Now, the generator of
$H_2(\pi_1(S),\Z)$
can be realized by an explicit cycle in the standard complex [B], that is,
$\sum^{2g}_{j=1} (I_{j-1}|x_j) + (I_{j-1}x_j|y_j) -(I_{j-1}x_j y_jx^{-1}_j|y_j) - (I_j|y_j)$
where
$I_j = [x_1,y_1] \ldots [x_j, y_j]$.
On the other hand, the universal Euler class may be realized by a cocycle
$A,B \mapsto \ell (A,B)$
as defined in Section 1. So the Euler number will be

$$\sum^y_{j=1} \ell (\bar I_{j-1}, \bar x_j) + \ell (\bar I_{j-1}, \bar x_j, \bar y_i) - \ell (\bar I_{j-1} \bar x_j \bar y_j \bar x^{-1}_j, \bar y_j) - \ell (\bar I_j, \bar y_j),$$
where
$\bar I_j$
is defined in an obvious way.
\enddemo

\demo{\bf{3. Proof of the Main Theorem}}
Consider a closed surface $S$ of genus 
$2n$. 
A map
$x_i \rightarrow a_i, y_i \rightarrow b_i, 1 \le i \le 2n$
defines a homomorphism from
$\pi_1(S)$
to
$SL_2(\R)$.
Indeed,
$\prod^n_{i=1} [a_i,b_i] \cdot \prod^{2n}_{i=n+1} [a_i,b_i] = h \cdot r^{-1} h^{-1}r =1$.
Next, the Euler number of this representation is computed as above, so
$f(a,b)$
is always an integer. Moreover, if
$f(a,b) = 4n-2$,
then the representation above is discrete and faithful.
\hfill Q.E.D.
\enddemo

\demo{\bf{$SU(1,n)$-case, I}}
Consider a standard action of
$SU(1,n)$
on the unit ball
$B \subset \C^n$
with the complex hyperbolic metric. Let
$\omega$
be the K\"ahler form of $B$. Fix a point
$\infty$
in the sphere at infinity and consider a function
$\varphi(A,B) = \psi(\infty, A(\infty), AB(\infty))$,
where
$\psi(x,y,z)$
is an integral of
$\omega$
over any surface, spanning the geodesic triangle with vertices
$x,y,z$.
Let
$a_i,b_i, \, 1\le i \le g$
be matrices in
$SU(1,n)$.
Let
$h = \prod^g_{i=1} [a_i, b_i]$
and suppose $h$ has a nonisotropic eigenvector. Then there exists a reflection $r$, commuting with $h$. Define
$a_i, b_i, \, i \ge g+1$
as in Introduction.
\enddemo

\proclaim{3.1 Theorem}
Define
$f(a,b)$
by the formula in the Main Theorem. Then
\roster
\item"(a)" $f(a,b)$
is an integer and
$|f(a,b)| \le 2g-1$
\item"(b)" if
$f(a,b) = 2g-1$,
then
$\{a_i,b_i\}$
generate a discrete group in
$SU(1,n)$.
\endroster
\endproclaim

\demo{Proof}
Same as above with Toledo's theorem [To], [Re1] instead of Goldman's.
\enddemo

\demo{\bf {4. Bounded cohomology}}
The bounded cohomology theory is an invention of Mikhail Gromov. The idea is as follows: in the standard complex, computing the real group cohomology we look only at bounded cochains, that is, bounded functions
$f : \undersetbrace i \to {G \times G \times \ldots \times G} \rightarrow \R$.
The resulted cohomology spaces are called
$H^i_b(G,\R)$.
There is a canonical homomorphism
$H^i_b(G,\R) \rightarrow H^i(G,\R)$.
\enddemo

\demo{4.1 Example}
Let $M$ be a symmetric space of negative curvature with isometry group $G$. Let
$\omega$
be any $G$-invariant $i$-form on $M$. Then one gets a (Borel) class
$Bor \, (\omega) \in H^i_{cont}(G,\R)$
(see [Re 1] for example), which may be represented by bounded cocycle [Gr 1]. The Euler class in
$SU(1,n), \, n \ge 1$
is a further specialization.
\enddemo

\demo{4.2 Second bounded cohomology and combinatorial group theory}
Consider a kernel of the map
$H^2_b(G,\R) \rightarrow H^2(G,\R)$.
It gives rise to a function
$f:G \rightarrow \R$
satisfying
$|f(xy) -f(x)-f(y)| \le C$.
Moreover, this function may be chosen a class function, that is,
$f(xyx^{-1}) = f(y)$
and such that
$f(x^n) = nf(x)$
[BG].
\enddemo

Next, for an element
$z \in G^{\prime} = [G,G]$
a {\it genus norm} is the smallest integer $g$ such that $z$ is a product of $g$ generators. A theorem of Culler [C] states:

\proclaim{Theorem 4.2 (Culler)}
If $G$ is a f.g. free group.There is a positive constant such that for any
$z \in G^{\prime}$,
$\parallel z^n\parallel_{genus} \ge const \, \cdot n$.
\endproclaim

The following result has been proved by Gromov [Gr2] and author [Re3].

\proclaim{Theorem 4.3}
Let $G$ be geometrically hyperbolic, that is a fundamental group of a manifold of pinched negative curvature with
$i(x) \underset{x\rightarrow \infty} \to \longrightarrow \infty$
(e.g. compact). Then the conclusion of the Theorem 4.2 holds.Here $i(x)$ is the injectivity radius at the point $x$.
\endproclaim

Now, we have

\proclaim{Proposition 4.4}
Let
$f:G \rightarrow \R$
be as above. If
$f(z) \ne 0$,
then
$\parallel z^n \parallel_{genus} \ge const \, \cdot n$.
\endproclaim

\demo{Proof}
Let
$z^n = \prod^g_{i=1} [x_i,y_i]$.
Then

$$\gather |n \cdot f(z)| = |f(z^n)| = |f(\prod^g_{i=1} [x_i,y_i])| = |\sum^g_{i=1} f([x_i,y_i])| + \\
+ C \cdot g \le 3c \cdot g + C \cdot g = 4C \cdot g, \quad \text{so} \quad g \ge \frac{n|f(z)|} {4C}, \endgather$$
\hfill Q.E.D.
\enddemo

\demo{\bf{Genus norm in lattices in $SU(1,n)$}} \enddemo

\proclaim{Theorem (4.5)}
Let
$\Gamma \subset SU(1,n)$
be a lattice with
$H_2(\Gamma,\R) =0$. 
There exists a nonzero function
$f : \Gamma \rightarrow \R$
as above such that if
$f(z) \ne 0$,
then
$\parallel z^n \parallel_{genus} \ge const \, \cdot n$.
\endproclaim

\demo{Remarks}
1. Observe that if
$f(z) \ne 0$,
then for any
$y \in G$
and
$\kappa$
big enough,
$f(z^{\kappa}y) \ne 0$,
so there are ``many'' $z$ for which the Theorem 6.1 applies. \newline
2. If
$\Gamma$
is cocompact or if
$\underset{\gamma \rightarrow \infty} \to {|Tr \, \gamma|} \rightarrow \infty$
in
$\Gamma$
then the conclusion of the Theorem follows from 4.3.
\enddemo

The proof of theorem 4.5 will be completed in 5.2.

\demo{\bf{5. Ergodic cocycle and measurable transfer}}
Let $G$ be a locally compact group, $H$ a closed subgroup and
$X = G/H$.
Suppose $G$ has an invariant finite Borel measure
$\mu$
on $X$. Suppose we have a measurable section
$\delta : X \rightarrow G$.
For any
$g \in G$
and
$x \in X$
we define
$\lambda (g,x) \in H$
as unique element such that

$$s(gx) = gs(x) \lambda (g,x)$$
This defines a map of groups

$$G \overset {\lambda}\to \longrightarrow H^X$$
Now, $G$ acts on
$H^X$
by changing the argument. The map
$\lambda$
is well-known to be a {\it cocycle} for first non-abelian cohomology. Suppose
$f(h_1, \ldots, h_n)$
is a measurable $n$-cocycle on $H$. Then the composition
$f \circ \lambda : G^{(n)} \rightarrow \R$
is a measurable cocycle valued in the $G$-module of measurable functions on $X$. If $f$ is bounded, so is
$f \circ \lambda$.
In this case,
$\int_X f \cdot \lambda$
is a real bounded $n$-cocycle on $G$ and we have a well-defined map

$$H^n_b(H,\R) \overset t \to \longrightarrow H^n_b (G,\R)$$
Moreover, the composition
$H^n_b(G,\R) \overset {res} \to \longrightarrow H^n_b(H,\R) \overset t \to \longrightarrow H^n_b (G,\R)$
is a multiplication by
$\text{Vol} \, X$.
\enddemo

The proof of all this facts is easily adopted from [Gr1].See [Re4].

\demo{5.2}
As an immediate corollary, we state:
\enddemo

\proclaim{Theorem (5.2)}
Let
$\Gamma$
be a lattice in either 1)
$SO(n,1)$
or 2)
$SU(n,1)$
or 3)
$SU_{\H}(n,1)$.
Then in case 1)
$H^n_b(\Gamma) \ne 0$;
in case 2)
$H^{2\kappa}_b(\Gamma, \R) \ne 0$
for all
$1 \le \kappa \le n$;
in case 3)
$H^{4\kappa}_b(\Gamma,\R)$
for all
$1 \le \kappa \le n$.
\endproclaim

\demo{Proof}
Let us prove 2), since the rest is similar.
$SU(n,1)$
acts isometrically on the complex ball
$B^n$.
The K\"ahler form 
$\omega$
defined, as in 4.1, a class
$\omega$
in
$H^2_{b,cont}(SU(n,1),\R)$.
It is nontrivial since for any cocompact
$\Gamma$
the restriction on
$H^2(\Gamma)$
gives the K\"ahler class of
$B/\Gamma$.
Now for any
$\Gamma$,
the restriction on
$H^2_b(\Gamma,\R)$
must be nontrivial, other- wise
$\text{Vol} (SU(n,1)/\Gamma) \cdot \omega =0$
even as a class in
$H^2(SU^{\delta} (n,1),\R)$.
\enddemo

\demo{Proof of the theorem 4.5}
We need only to handle the case
$H_2(\Gamma,\R) =0$.
Then the restriction of the just defined class in
$H^2_b(\Gamma,\R)$
on
$H^2(\Gamma, \R)$
is zero, so by 4.2 we have a function $f$ with desired properties.
\enddemo


\enddocument